\input amstex
\documentstyle{amsppt}
\topmatter
\title
Quantum super spheres and their transformation groups, representations and little $t$-Jacobi polynomials
\endtitle
\author
Yi Ming Zou
\endauthor
\address
Department of Mathematical Sciences,
University of Wisconsin-Milwaukee,
Milwaukee, WI 53201, USA
\endaddress
\email
ymzou\@uwm.edu
\endemail
\abstract
Quantum super 2-shpheres and the corresponding quantum super transformation group are introduced in analogy to the well-known quantum 2-shpheres and quantum $SL(2)$, connection between little $t$-Jacobi polynomials and the finite dimensional representations of the quantum super group is formulated, and the Peter-Weyl theorem is obtained.
\endabstract
\endtopmatter

\document
The quantum group $SU_{q}(2)$ introduced by Woronowicz and the quantum 2-spheres introduced by Podle\`s play basic roles in quantum group and quantum homogeneous space theory.  In terms of ``quantum homogeneous spaces", the Podle\`s spheres can be regarded as quantum homogeneous spaces with $SU_{q}(2)$ as their transformation group.  The quantum group $SU_{q}(2)$ (or $SL_{q}(2)$) and quantum 2-spheres have been studied intensively (there is a comprehensive bibliography in [KS]). The results obtained in these studies revealed connections with other fields of mathematics and form the foundation of more general studies in quantum group and quantum homogeneous spaces.  
\par
For a contragredient Lie superalgebra, besides $sl(2)$, there are two other basic building blocks, namely $sl(1,1)$ and $osp(1,2)$.  The structure of the corresponding quantized enveloping algebra of $sl(1,1)$ is relatively simple because of the condition $f^{2}=e^{2}=0$ on the generators. Though there have been many publications devoted to the quantized enveloping algebra $U_{q}(osp(1,2n))$ of $osp(1,2n)$ (especially when $n=1$) and its representations, there does not seem to be a detail study of the corresponding quantum group that is comparable to the existing $SU_{q}(2)$ theory, perhaps due partially to the fact that in general, the theory of $osp(1,2n)$ is similar to that of $so(2n+1)$, and for the quantum $osp(1,2)$ the theory is similar to that of $SU_{q}(2)$. But the similarity deserves further explanation. Over the field of complex numbers, the first nontrivial representation of the Lie superalgebra $osp(1,2)$ is 3-dimensional, and all nontrivial finite dimensional representations have odd dimensions. In contrary, the finite dimensional representations of the deformation $U_{q}(osp(1,2))$ assume both even and odd dimensions (see [MZ]).  This makes the representation theory of the $U_{q}(osp(1,2))$ even closer to that of $U_{q}(sl(2))$ in this aspect.  On the other hand, there is a difference between the nontrivial even dimensional representations of these two quantum algebras: the ones for $U_{q}(sl(2))$ are self-dual, but not for $U_{q}(osp(1,2))$.  When the duals of these modules are used to construct the quantum groups, the element corresponds to the quantum determinant is not in the center in the $U_{q}(osp(1,2))$ case (see section 1 below). Nevertheless, a parallel theory can still be developed, which is the purpose of this paper.  Since the theory of $SL_{q}(2)$ is treated (more or less) in any book on quantum groups, we shall present our discussion in a standard way without pointing out the original sources of the similar results in the non-super case.  For references of the quantum 2-spheres, besides Podle\`s original paper [Po], we also refer the readers to references [DK], [Ka] and [KS], and we shall use notation that similar to those used in the later publications.  It should be pointed out that the quantum super spheres discussed in this paper do not contain the Podle\`s quantum spheres as a special case, instead, they are parallel to the Podle\`s spheres in the super setting. The approach to the finite dimensional comodules of the quantum super group (the transformation group of the super spheres), the matrix elements of these comodules and their connection with the little $t$-Jacobi polynomials, as well as the Peter-Weyl theorem, are due to [MNU]. 
\par
One can see why such a parallel theory holds by means of the Radford-Majid bosonization (The author would like to thank the referee for suggesting this explanation.).  Start with the coordinate algebra of $2\times 2$ quantum matrices $\Cal{O}_{t}(M_{2}(\Bbb{C})) = \Bbb{C}<A,B,C,D>$ with the relations
$$
\gathered
AB=tBA, \quad AC=tCA, \quad BD=tDB, \quad CD=tDC, \quad CB=BC\\
AD-DA=(t-t^{-1})BC.
\endgathered
$$
Let $\Bbb{Z}_{2}= <g>$ act as automorphisms of $\Cal{O}_{t}(M_{2}(\Bbb{C}))$ by
$$
g\cdot \left( \matrix A & B\\ C & D \endmatrix \right)=
\left( \matrix A & -B\\ -C & D \endmatrix \right).
$$
Then $\Cal{O}_{t}(M_{2}(\Bbb{C}))\#\Bbb{C}\Bbb{Z}_{2}$ becomes a bialgebra with coalgebra structure as \newline $\Cal{O}_{t}(M_{2}(\Bbb{C}))\otimes\Bbb{C}\Bbb{Z}_{2}$.  There is a bialgebra surjection
$$
\pi : \Cal{O}_{t}(M_{2}(\Bbb{C}))\#\Bbb{C}\Bbb{Z}_{2}\longrightarrow \Bbb{C}\Bbb{Z}_{2}
$$
sending $\left( \matrix A & B\\ C & D \endmatrix \right)$ to $\left( \matrix 1 & 0\\ 0 & g \endmatrix \right)$ and which is identity on $\Bbb{C}\Bbb{Z}_{2}$.  Let $R$ be the right $\Bbb{C}\Bbb{Z}_{2}$-coinvariants of $\Cal{O}_{t}(M_{2}(\Bbb{C}))\#\Bbb{C}\Bbb{Z}_{2}$.  By bossonization, it has the structure of a braided bialgebra in the category of left $\Bbb{C}\Bbb{Z}_{2}$ Yetter-Drinfeld modules.  If we put
$$
a=A\#1,\quad b= B\#g, \quad c=C\#1, \quad d=D\#g,
$$
we have $R=\Bbb{C}<a,b,c,d>$ and the relations (1.1) and (1.2) (see section 1 below) hold, with the coproduct of $R$ defined so that (see [Mo, 10.6])
$$
\Delta(b\#h)=\sum b^{1}\#(b^{2})_{-1}h_{1}\otimes (b^{2})_{0}\#h_{2}.
$$
Further $a,d$ have degree $0$ and $b,c$ have degree 1.  The action of $\Bbb{Z}_{2}$ is given by $g\cdot a = a$, $g\cdot b = -b$, $g\cdot c = -c$, and $g\cdot d = d$.  Thus the braided bialgebra $R$ coincides with the $\Bbb{Z}_{2}$-graded algebra $\Cal B$ (see section 1).  Now $\delta = AD -tBC$ is a central group-like element in $\Cal{O}_{t}(M_{2}(\Bbb{C}))\#\Bbb{C}\Bbb{Z}_{2}$.  If we put $\sigma = \delta g$, then we have
$$
\Cal{O}_{t}(M_{2}(\Bbb{C}))/(\delta^{2}-1)\#\Bbb{C}\Bbb{Z}_{2} = \Cal{A}(\sigma)\#\Bbb{C}\Bbb{Z}_{2},
$$
i.e. $\Cal{A}(\sigma)$ is the braided Hopf algebra corresponding to the above quotient Hopf algebra.
\par
There is a different discussion of a certain type of quantum supergroup associated with $osp(1,2n)$ in [LZ]. There the quantum supergroup is the dual to the so-called type 1 deformations of the simple modules of $osp(1,2n)$. In the $osp(1,2)$ case, as mentioned before, these representations assume only odd dimensions.
\par
In section 1, we define the $\Bbb{Z}_{2}$-graded algebras $\Cal B$ and $\Cal{A}(\sigma)$, they can be considered as the $\Bbb{Z}_{2}$-graded coordinate ring and the corresponding quantum super group ($\Bbb{Z}_{2}$-graded Hopf algebra) of a certain type of quantum super $2\times 2$ matrices (see section 2 and section 3).  In order to obtain the quantum super group from the coordinate ring $\Cal B$, we need to invert a group-like element which is not in the center, this explains the existence of the element $\sigma$.  In section 3, we show that these $\Bbb{Z}_{2}$-graded algebras arise naturally from the transformations of certain quantum super space. In section 4, we discuss the duality between $\Cal B$ (resp. $\Cal{A}(\sigma)$) and $U_{q}(osp(1,2))$, the dual pairing is provided by a 2-dimensional $U_{q}(osp(1,2))$-module which requires the use of $t=\bold{i}\sqrt{q}$.  Hence though the deformation parameter of the quantized enveloping algebra is $q$, the deformation parameter for the quantum super group is $t$.  In section 4, we define some quantum super 2-spheres, they are constructed as the quantum super homogeneous spaces of $\Cal{A}(\sigma)$.  The construction is similar to that of [Po].  In section 5, we study the finite dimensional comodules of $\Cal{A}(\sigma)$.  The description of all irreducible finite dimensional comodules is given in Theorem 5.1, which also gives the decomposition of $\Cal{A}(\sigma)$ as direct sum of irreducible bi-subcomdules.  In section 6, we provide formulas for the matrix coefficients of the irreducible comodules in terms of little $t$-Jacobi polynomials.  The powers of $\bold{i}=\sqrt{-1}$ in the formulas are coming from the square roots of some powers of $-1$, which signify the super-commutativity. In section 7, we derive the Peter-Weyl theorem for $\Cal{A}(\sigma)$.
\par
\head
1. The algebras $\Cal{B}$ and ${\Cal A}(\sigma)$
\endhead
\par
Let ${\Bbb C}$ be the complex number field and let $q\in \Bbb{C}$.  We assume through out that $q\ne 0$ and $q$ is not a root of $1$.  Let $\bold{i} = \sqrt{-1}$, and let $t =\bold{i}\sqrt{q}$. We define a $\Bbb{Z}_{2}$-graded associative algebra $\Cal B$ over $\Bbb C$ with generators $a, b, c, d$ (and $1$) subject to the relations
$$
\gathered
ab = tba, \qquad ac = tca, \qquad bc = -cb, \\
bd = -tdb, \qquad cd = -tdc, \qquad ad - da = (t^{-1}-t)bc.
\endgathered \tag 1.1
$$
The parities of the elements of $\Cal B$ is given by $p(a)=p(d)=0$ and $p(b)=p(c)=1$, where $p$ is the parity function. We define a 
$\Bbb{Z}_{2}$-graded coalgebra structure on $\Cal B$ by 
$$
\gathered
\Delta \left( \matrix a & b\\ c & d \endmatrix \right)=
\left( \matrix a & b\\ c & d \endmatrix \right)\otimes
\left( \matrix a & b\\ c & d \endmatrix \right),\\
\varepsilon (a)=\varepsilon(d) = 1, \quad \varepsilon(b)=\varepsilon(c)=0.
\endgathered \tag 1.2
$$
\par
It is easy to check that the element $ad + tbc$ is a group-like element in $\Cal B$, it commutes with $a$ and $d$, but 
$$
 b(ad + tbc)= -(ad + tbc)b,\qquad c(ad + tbc)= -(ad + tbc)c. 
$$
\par
Thus we introduce another algebra $\Cal{B}(\sigma)$ by adding an even generator $\sigma$ to those of $\Cal B$ and the following extra relations
$$
\sigma a = a\sigma, \quad \sigma d = d\sigma, \quad
\sigma b = -b\sigma, \quad \sigma c = -c\sigma, \quad \sigma^{2}=1, \tag 1.3
$$
$$
\Delta(\sigma)=\sigma\otimes\sigma, \qquad \varepsilon(\sigma)=1. \tag 1.4
$$
It is straightforward to check that the algebra $\Cal{B}(\sigma)$ is a $\Bbb{Z}_{2}$-graded bialgebra with generators $a,b,c,d,\sigma$ and generating relations (1.1)-(1.4). 
\par
We now let $\Cal{A}(\sigma) = \Cal{B}(\sigma)/(ad + tbc = \sigma)$, then $\Cal{A}(\sigma)$ is also a $\Bbb{Z}_{2}$-graded bialgebra.  We define a $\Bbb{Z}_{2}$-graded anti-automorphism $S$ of $\Cal{A}(\sigma)$ by
$$
S\left( \matrix a & b\\ c & d \endmatrix \right)=
\left( \matrix d\sigma & -t^{-1}b\sigma\\ tc\sigma & a\sigma \endmatrix \right), \qquad 
S(\sigma)= \sigma, \tag 1.4
$$
$$
S(xy) = (-1)^{p(x)p(y)}S(y)S(x), \quad \text{$x,y \in \Cal{A}(\sigma)$}. \tag 1.5
$$
It is easy to see that $S$ preserves the generating relations (1.1), (1.3) and $ ad + tbc = \sigma$, therefore $S$ is well-defined. Let $\Bbb{Z}_{+}=\{0,1,2, ...\}$. 
\par
\proclaim{Theorem 1.1}
With the anti-automorphism $S$ defined by (1.4) and (1.5), the bialgebra $\Cal{A}(\sigma)$ is a $\Bbb{Z}_{2}$-graded Hopf algebra over $\Bbb C$. A basis of $\Cal{A}(\sigma)$ is provided by the monomials $a^{i}b^{j}c^{k}d^{l}\sigma^{s}$, $i,j,k,l \in \Bbb{Z}_{+}$, $i=0$ or $l=0$, and $s=0,1$. 
\endproclaim
\demo{Proof}
For the first statement, we only need to verify the antipode property for $S$ on the generators and this is clear from the formulas. For the second statement, we apply the diamond lemma. 
\enddemo
\par
If $q < 0$ is a real number, then $t$ is a real number, and we can define a $\ast$-algebra structure on $\Cal{A}(\sigma)$ by 
$$
\ast \left( \matrix a & b\\ c & d \endmatrix \right)=
\left( \matrix d\sigma & tc\sigma\\ -t^{-1}b\sigma & a\sigma \endmatrix \right),\qquad \ast(\sigma)=\sigma, \tag 1.6
$$
together with the requirement that $\ast$ be an anti-multiplicative and anti-linear map.
\par
\proclaim
{Lemma 1.2} Under the assumption that $q < 0$, we have
$$\gather
\ast^{2} = id, \quad (\ast\otimes\ast)\circ \Delta = 
\Delta\circ\ast, \quad \varepsilon(\ast(x))=\overline{\varepsilon(x)} \quad\text{($x\in \Cal{A}(\sigma)$)}, \tag 1.7 \\
\ast\circ S\circ\ast\circ S = S\circ\ast\circ S \circ\ast = id. \tag 1.8
\endgather
$$
\endproclaim
\demo{Proof}
The proof is given by straightforward computations.
\enddemo
\par
By this lemma, if $q < 0$, then $\Cal{A}(\sigma)$ is a $\Bbb{Z}_{2}$-graded Hopf $\ast$-algebra, and we can define its compact real form as in the non-super case. 
\par
Let $\Cal{O}(K)$ denote the Hopf algebra $\Bbb{C}[z,z^{-1}]$ with comultiplication $\Delta(z)=z\otimes z$, counit $\varepsilon(z)=1$ and antipode $S(z)=z^{-1}$. Define a Hopf algebra homomorphism $\phi_{K}: \Cal{A}(\sigma)
\longrightarrow \Cal{O}(K)$ by
$$
\phi_{K}(a)=z, \quad \phi_{K}(b)=\phi_{K}(c)=0, \quad \phi_{K}(d)=z^{-1},
\quad \phi_{K}(\sigma)=1. \tag 1.9
$$
Define a left $\Cal{O}(K)$-comodule structure $L_{K}$ and a right 
$\Cal{O}(K)$-comodule structure $R_{K}$ on $\Cal{A}(\sigma)$ by
$$
L_{K}=(\phi_{K}\otimes id_{\Cal{A}(\sigma)})\circ\Delta, \qquad
R_{K}=(id_{\Cal{A}(\sigma)}\otimes \phi_{K})\circ\Delta. \tag 1.10
$$
For $m,n\in \Bbb{Z}$ set
$$
\Cal{A}(\sigma)[m,n] = \{x\in\Cal{A}(\sigma) : L_{K}(x)=z^{m}\otimes x
\quad\text{and} \quad R_{K}(x)=x\otimes z^{n}\}. \tag 1.11
$$
It is clear that $\Cal{A}(\sigma)[m,n] \cdot \Cal{A}(\sigma)[r,s]\subset
\Cal{A}(\sigma)[m+r,n+s]$, hence $\Cal{A}(\sigma)[0,0]$ is a $\Bbb{C}$-subalgebra of $\Cal{A}(\sigma)$ and each $\Cal{A}(\sigma)[m,n]$ is an $\Cal{A}(\sigma)[0,0]$-bimodule.  For each pair $(m,n)\in \Bbb{Z}^{2}$ with $m\equiv n (\mod 2)$, define $e_{mn}\in \Cal{A}(\sigma)[m,n]$ by
$$
\gathered
e_{mn}=a^{(m+n)/2}c^{(n-m)/2},\qquad \text{if $m+n\ge 0, m\le n$},\\
e_{mn}=a^{(m+n)/2}b^{(m-n)/2},\qquad \text{if $m+n\ge 0, m\ge n$},\\
e_{mn}=b^{(m-n)/2}d^{(-n-m)/2},\qquad \text{if $m+n\le 0, m\ge n$},\\
e_{mn}=c^{(n-m)/2}d^{(-n-m)/2},\qquad \text{if $m+n\le 0, m\le n$}.
\endgathered\tag 1.12
$$
\par
By Theorem 1.1, the following proposition holds (cf. [MNU, 1.2]).
\proclaim{Proposition 1.3} (1) $\Cal{A}(\sigma)[0,0]$ is the subalgebra  $\Bbb{C}<bc,\sigma>$ of $\Cal{A}(\sigma)$ generated by $bc$ and $\sigma$.
\par
(2) $\Cal{A}(\sigma) = \oplus_{m,n\in\Bbb{Z}}\Cal{A}(\sigma)[m,n]$.
\par
(3) $\Cal{A}(\sigma)[m,n] \ne 0$ if and only if $m\equiv n (\mod 2)$, and if $m\equiv n (\mod 2)$, $\Cal{A}(\sigma)[m,n]$ is a free left (or right) 
$\Cal{A}(\sigma)[0,0]$-module of rank one with basis $e_{mn}$.
\endproclaim
\par
\head 
2. Geometric realization of $\Cal B$
\endhead
\par
We consider a geometric approach to the algebra $\Cal B$.  Let $\Cal{O}(\Bbb{C}^{0,1}_{t})$ be the $\Bbb{Z}_{2}$-graded algebra generated by $x$ and $y$ satisfying the relation $xy=tyx$, let the parities be given by $p(x)=0, p(y)=1$.
\par
Let $(\Cal{B'}, \Delta, \varepsilon)$ be a $\Bbb{Z}_{2}$-graded bialgebra over $\Bbb C$ with generators $a',b',c',d'$ and the parities defined by $p(a')=p(d')=0, p(b')=p(c')= 1$ (these parity assignments need to be met in order to define the actions in (2.1)).  We define the left and the right actions of $\Cal B'$ on the vectors $(x, y)^{T}$ and $(x, y)$ by 
$$\gathered
\left( \matrix a' & b'\\ c' & d' \endmatrix \right)\otimes
\left(\matrix x \\ y \endmatrix \right)=
\left( \matrix a'\otimes x + b'\otimes y\\ c'\otimes x + d'\otimes y \endmatrix \right),\\
\left(\matrix x, y\endmatrix \right)\otimes
\left( \matrix a' & b'\\ c' & d' \endmatrix \right)=
\left( \matrix x\otimes a' + y\otimes c', x\otimes b' + y\otimes d' \endmatrix \right).
\endgathered\tag 2.1 
$$
\par
Recall that if $(A, \Delta, \varepsilon)$ is a $\Bbb{Z}_{2}$-graded bialgebra, then a quantum $A$-space (left space, right space can be defined similarly) is given by a pair $(B, \psi)$, where $B$ is an 
associative $\Bbb{Z}_{2}$-graded algebra and $\psi$ is an (even) algebra homomorphism $B\longrightarrow A\otimes B$ such that 
$$
(id\otimes \psi)\circ\psi = (\Delta\otimes id)\circ\psi, \qquad
(\varepsilon\otimes id)\circ\psi = id. \tag 2.2
$$
\par
\proclaim{Theorem 2.1}
The actions defined by (2.1) provide both a left and a right (respectively) quantum $B'$-space structures for $\Cal{O}(\Bbb{C}^{0,1}_{t})$ if and only if $a',b',c',d', \Delta, \varepsilon$ satisfy (1.1) and (1.2).
\endproclaim
\demo{Proof}
 Note that for the left action, the second equality in (2.2) is equivalent to $\varepsilon(a')x + \varepsilon(b')y = x$ and 
$\varepsilon(c')x + \varepsilon(d')y = y$ under the assumption that $\psi$ is an algebra homomorphism, and these equalities are equivalent to the conditions on $\varepsilon$ in (1.2) by the linearly independence of $x$ and $y$.  The same is also true for the right action.  Under the assumption that $\psi$ is an algebra homomorphism, the first condition in (2.2) is equivalent to 
$$
(id\otimes \psi)\circ\psi (x) = (\Delta\otimes id)\circ\psi (x), \quad
(id\otimes \psi)\circ\psi (y)= (\Delta\otimes id)\circ\psi (y),
$$
i.e.
$$\gather
(a'\otimes a' + b'\otimes c')\otimes x + (a'\otimes b' + b'\otimes d')\otimes y
= \Delta (a')\otimes x + \Delta(b')\otimes y,\\
(c'\otimes a' + d'\otimes c')\otimes x + (c'\otimes b' + d'\otimes d')\otimes y
= \Delta (c')\otimes x + \Delta(d')\otimes y.
\endgather
$$
These conditions are equivalent to the conditions for $\Delta$ in (1.2).
\par
It remains to prove that $\psi$ is a $\Bbb{Z}_{2}$-algebra homomorphsim is equivalent to the conditions in (1.1).  Consider the left action.  In order for $\psi$ to be an algebra homomorphsim, it is necessary and sufficient that 
$$
(a'\otimes x + b'\otimes y)(c'\otimes x + d'\otimes y) =
t(c'\otimes x + d'\otimes y) (a'\otimes x + b'\otimes y),
$$
i.e. (use the given parities on $a',b',c',d',x,y$ and $xy=tyx$)
$$
(a'c'-tc'a')\otimes x^{2} + (a'd'-tc'b'-t^{-1}b'c'-d'a')\otimes xy +
(b'd'+td'b')\otimes y^{2} = 0.
$$
Which is equivalent to 
$$
a'c'=tc'a', \quad a'd'-tc'b' = d'a'+t^{-1}b'c',\quad b'd'=-td'b'. \tag 2.3
$$
Similarly, the right action leads to  
$$
a'b'=tb'a', \quad a'd'+tb'c'=d'a'-t^{-1}c'b',\quad c'd'=-td'c'. \tag 2.4
$$
The two middle equalities imply $t(b'c'+c'b')=-t^{-1}(c'b'+b'c')$, i.e.
$(t+t^{-1})(b'c'+c'b')=0$, or $b'c'=-c'b'$ by our definition of $t$.  The rest of the proof is clear now.
\enddemo
\par
By the above theorem, the quantum superspace $\Cal{O}(\Bbb{C}^{0,1}_{t})$ is a $\Cal B$-space.  Furthermore, the proof of the theorem shows that 
$\Cal{O}(\Bbb{C}^{0,1}_{t})$ is both a left and a right comodule algebra of $\Cal B$.  If we add an extra relation $y^{2}=0$ to the definition of $\Cal{O}(\Bbb{C}^{0,1}_{t})$, the resulting quantum superspace $\Cal{O}(\Bbb{C}^{0,1}_{t})/(y^{2})$ is one of the quantum superspaces defined in [Ma, p.136]. However, this later space is not a $\Cal B$-space with the actions defined by (2.1).
\par
\head
3. The dual of $\Cal B$
\endhead
\par
Recall that the natural $\Bbb{Z}_{2}$-graded algebra structure on the 
dual space $\Cal{B}^{\ast} = \text{Hom}_{\Bbb{C}}(\Cal{B}, \Bbb{C})$ is defined by 
$$
(fg)(x)=(f\otimes g)\circ\Delta (x), \qquad f,g\in \Cal{B}^{\ast}, x\in \Cal{B},
\tag 3.1
$$
and
$$
(f\otimes g)(x\otimes y) = (-1)^{p(x)p(g)}f(x)g(y), \qquad f,g\in \Cal{B}^{\ast}, x,y\in \Cal{B}. \tag 3.2
$$
\par
Let $k^{\pm 1}, e, f \in \Cal{B}^{\ast}$ be defined by
$$
\gathered
k^{\pm 1}\left( \matrix a & b\\ c & d \endmatrix \right)=
\left( \matrix t^{\pm 1} & 0\\ 0 & -t^{\mp 1} \endmatrix \right),\\
e\left( \matrix a & b\\ c & d \endmatrix \right)=
\left( \matrix 0 & \frac{t-t^{-1}}{q-q^{-1}}\\ 0 & 0 \endmatrix \right), \quad
f\left( \matrix a & b\\ c & d \endmatrix \right)=
\left( \matrix 0 & 0\\ 1 & 0 \endmatrix \right),
\endgathered\tag 3.3
$$
together with the requirements that $k^{\pm 1}$ are algebra homomorphisms and 
$$
\gathered
e(xy)= e(x)\varepsilon (y) + (-1)^{p(x)}k(x)e(y), \qquad e(1)=0,\\
f(xy)= f(x)k^{-1}(y) + (-1)^{p(x)}\varepsilon(x)f(y), \qquad f(1)=0,
\endgathered\tag 3.4
$$
for all $x,y\in \Cal{B}$.
\par
We let $\Cal U$ be the $\Bbb{Z}_{2}$-graded subalgebra of $\Cal{B}^{\ast}$ generated by $k^{\pm 1}, e,f$. Note that the $\Bbb{Z}_{2}$-grading of $\Cal U$ is given by $p(k^{\pm 1})=0$ and $p(e)=p(f)=1$. For our convenience, we let $X = \{a,b,c,d \}$.
\par
\proclaim{Lemma 3.1}
The generators of $\Cal U$ satisfy the relations
$$\gathered
kk^{-1}=k^{-1}k=1,\\
kek^{-1}=qe, \qquad kfk^{-1}=q^{-1}f,\\
ef + fe = \frac{k-k^{-1}}{q-q^{-1}}.
\endgathered\tag 3.5
$$
\endproclaim
\demo{Proof}
Note that elements of the form $x = x_{1}\cdots x_{m}$, $x_{i}\in X$ ($1\le i \le m$), span $\Cal B$ as a vector space, so by linearity we only need to verify the identities on these elements.  It is clear that if $x\in X$, then the identities hold.
\par
To prove $ kk^{-1}=k^{-1}k=1$, we observe that if one of the $x_{i}$'s that appear in  $x = x_{1}\cdots x_{m}$ is $b$ or $c$, then by the definition of $k^{\pm 1}$, we have $k^{\pm 1}(x)=0$.  Therefore by (1.2), if $b$ or $c$ appears in $x$, we also have $kk^{-1}(x)=0=\varepsilon(x)$.  Thus we can assume that the $x_{i}$'s are either $a$ or $d$.  Then by (1.2) we have 
$$\split
kk^{-1}(x)&=kk^{-1}(x_{1}\cdots x_{m})= 
k(x_{1}\cdots x_{m})k^{-1}( x_{1}\cdots x_{m})\quad\text{(other terms are 0)}\\
&=k(x_{1})\cdots k(x_{m}) k^{-1}(x_{1})\cdots k^{-1}(x_{m})
= k(x_{1})k^{-1}(x_{1})\cdots k(x_{m})k^{-1}(x_{m})\\
&=\varepsilon(x_{1})\cdots\varepsilon(x_{m})=\varepsilon(x).
\endsplit
$$
Hence $kk^{-1}=1$. Simlarly, $k^{-1}k=1$.
\par
To prove the other three identities, we proceed by induction on the length of 
$x$. Assume both sides of the identities are equal when evaluated at $x$, we verify that they are also equal when evaluated at $xy$ for $y\in X$.  We shall give detail on two cases, the other cases are similar.
\par
For $kek^{-1}=qe$, we verify that
$(kek^{-1})(xb)=qe(xb)$.  Let 
$$
\Delta(x)= \sum x_{1}\otimes x_{2},\qquad \Delta(x_{1})=\sum x_{11}\otimes x_{12}
$$
as usual.  We have
$$
\split
kek^{-1}(xb) &= (ke\otimes k^{-1})(\sum(x_{1}\otimes x_{2})(a\otimes b + b\otimes d))\\
&= \sum (-1)^{p(x_{2})}(ke)(x_{1}b)k^{-1}(x_{2}d)\\
&= \sum (-1)^{p(x_{2})}(k\otimes e)((x_{11}\otimes x_{12})(a\otimes b + b\otimes d))k^{-1}(x_{2}d)\\
&= \sum (-1)^{p(x_{2})+ p(x_{11})}k(x_{11}a)e(x_{12}b)k^{-1}(x_{2}d)\\
&= \sum (-1)^{p(x_{2})+ p(x_{11}) + p(x_{12})}
k(x_{11}a)k(x_{12})e(b)k^{-1}(x_{2}d) \quad \text{(by (3.4))}\\
&= \sum(-1)^{p(x)}qk(x_{11})k(x_{12})e(b)k^{-1}(x_{2}) 
\quad\text{($k(a)k^{-1}(d)= q$)}\\
&= (-1)^{p(x)}qk(x)e(b)=qe(xb). \quad\text{(by (3.1) and (3.4))}
\endsplit
$$
\par
For $ef+fe= \frac{k-k^{-1}}{q-q^{-1}}$, we verify
$(ef+fe)(xb)= (\frac{k-k^{-1}}{q-q^{-1}})(xb)$ under the assumption that $kek^{-1}=qe$ and $kfk^{-1}=q^{-1}f$ have been proved.  
Since the right hand side is equal to 0, we only need to verify that the 
left hand side is also equal to 0. By (3.3) and (3.4), $(e\otimes f)(x_{1}a\otimes x_{2}b)=0$, thus
$$
\split
ef(xb) &= \sum (-1)^{p(x_{2})+p(x_{2}b)+p(x_{1})}
k(x_{1})e(b)f(x_{2})k^{-1}(d)\\
&= \sum (-1)^{p(x) + 1 + p(x_{1})}k(x_{1})f(x_{2})e(b)k^{-1}(d)\\
& = (-1)^{p(x)+1}(kf)(x)(ek^{-1})(b)\qquad\text{(since $ek^{-1}(b)=e(b)k^{-1}(d)$)}\\
&= -(kf\otimes ek^{-1})(x\otimes b). 
\endsplit
$$
Similarly, $fe(xb)=(fk\otimes k^{-1}e)(x\otimes b)$.  Since $kek^{-1}=qe$ and $kfk^{-1}= q^{-1}f$, we have $kf\otimes ek^{-1} = fk\otimes k^{-1}e$, and hence
$(ef+fe)(xb)=0$.
\enddemo
\par
Recall that the $q$-deformation $U_{q}(osp(1,2))$ of the universal enveloping algebra of the Lie superalgebra $osp(1,2)$ (see for example [MZ]) is generated by $K^{\pm 1}, E, F$ subject to generating relations as in (3.5).  Thus there is an onto algebra homomorphism $\varphi$ from $U_{q}(osp(1,2))$ to $\Cal U$.  Furthermore, there is a natural $\Bbb{Z}_{2}$-graded Hopf algebra structure on $\Cal U$ such that $\varphi$ is a $\Bbb{Z}_{2}$-graded Hopf algebra homomorphism.  The comultiplication, the antipode and counit on $\Cal U$ are provided by 
$$\gather
\Delta (k^{\pm 1}) = k^{\pm 1}\otimes k^{\pm 1}, 
\Delta (e) = e\otimes 1 + k\otimes e,
\Delta (f) = f\otimes k^{-1} + 1\otimes f; \tag 3.6\\
S(k^{\pm 1}) = k^{\mp 1}, \quad S(e)=-k^{-1}e, \quad S(f) = -fk; \tag 3.7\\
\varepsilon (k^{\pm 1})=1, \quad \varepsilon (e)=\varepsilon (f) = 0. \tag 3.8
\endgather
$$
\par
The coalgebra structure on $\Cal U$ is dual to the algebra structure of $\Cal B$, i.e. we have
$$
\Delta (f)(x\otimes y) = f(xy), \quad \varepsilon(f) = f(1), \quad f\in \Cal{U},
\quad x,y\in \Cal{B}. \tag 3.9
$$
\par
The homomorphism $\varphi$ induces a bilinear pairing between $U_{q}(osp(1,2))$ and $\Cal B$, this bilinear pairing satisfies (3.1)-(3.3) and (3.9). Using an argument based on the representation theory of $U_{q}(osp(1,2))$ (cf. [KS, p.122]), we can prove that this bilinear pairing is nondegenerate, and as a consequence, we conclude that $\varphi$ is injective.  Thus we have the following
\proclaim{Theorem 3.2}
The algebras $\Cal U$ and $U_{q}(osp(1,2))$ are isomorphic as $\Bbb{Z}_{2}$-graded Hopf algebras.
\endproclaim
\par
Note that we can define a pairing between $\Cal U$ and $\Cal{A}(\sigma)$ by setting 
$$
k^{\pm 1}(\sigma)=-1, \qquad e(\sigma)=f(\sigma)=0. \tag 3.10
$$
With this pairing, we have
$$
S(f)(x)=f(S(x)),\qquad f\in \Cal{U}, \quad x\in \Cal{A}(\sigma).\tag 3.11
$$
\par
\head
4. Quantum super spheres
\endhead
In this section, we define some super quantum homogeneous spaces of $\Cal{A}(\sigma)$, following [Po], we call them quantum super spheres. Let $M=(m_{ij})$ be the matrix defined by 
$$
M = \left(\matrix
a^{2} & \sqrt{1-q^{-1}}ab & \bold{i}b^{2}\\
\sqrt{1-q^{-1}}ac & ad + t^{-1}cb & \bold{i}\sqrt{1-q}db\\
\bold{i}c^{2} & -\bold{i}\sqrt{1-q}dc & d^{2}
\endmatrix\right). \tag 4.1
$$
\par
\proclaim{Lemma 4.1} We have
$$
\Delta(M) = M\otimes M, \quad \varepsilon(M) = (\delta_{ij}), \quad
S(M) = (\ast(M))^{T} \quad \text{if $q<0$}.  \tag 4.2
$$
That is, $M$ is a matrix corepresentation of $\Cal{A}(\sigma)$, and $M$ is unitary if $q<0$.
\endproclaim
\demo{Proof}
The proof is given by computations.  For example, by (1.1) and (1.2),
$$
\split
&\Delta(ad+t^{-1}cb)\\
&= (a\otimes a + b\otimes c)(c\otimes b + d\otimes d)
+ t^{-1}(c\otimes a + d\otimes c)(a\otimes b + b\otimes d)\\
&= ac\otimes ab + ad\otimes ad - bc\otimes cb + bd\otimes cd\\
&\qquad + t^{-1}(ca\otimes ab + cb\otimes ad + da\otimes cb - db\otimes cd)\\
&=(1+t^{-2})ac\otimes ab + (ad+t^{-1}cb)\otimes ad + (da+tcb)\otimes t^{-1}cb + (1+t^{2})db\otimes dc\\
&= (1-q^{-1})ac\otimes ab + (ad+t^{-1}cb)\otimes (ad+t^{-1}cb) + (1-q)db\otimes dc\\
&= (\sqrt{1-q^{-1}}ac, ad+t^{-1}cb, \bold{i}\sqrt{1-q}db)\otimes \left(\matrix
\sqrt{1-q^{-1}}ab\\
ad + t^{-1}cb\\
-\bold{i}\sqrt{1-q}dc \endmatrix\right).
\endsplit
$$
The others cases are similar.
\enddemo
\par
For $\alpha = (\alpha_{-1}, \alpha_{0}, \alpha_{1})\in \Bbb{C}^{3}\setminus(0)$, define $\bold{x}(\alpha) = (x(\alpha)_{-1}, x(\alpha)_{0}, x(\alpha)_{1}) = \alpha\cdot M$.  Abusing notation, we let $\bold{x}(\infty) = (x(\infty)_{-1}, x(\infty)_{0}, x(\infty)_{1})$ be defined by
$$
x(\infty)_{-1} = \left(\frac{t+t^{-1}}{t-t^{-1}}\right)^{\frac{1}{2}}ac,
 x(\infty)_{0} = ad + t^{-1}cb,  
x(\infty)_{1}=\left(\frac{t+t^{-1}}{t-t^{-1}}\right)^{\frac{1}{2}}db. \tag 4.2
$$
We denote by $S_{\alpha}$ the subalgebra of $\Cal{A}(\sigma)$ generated by $ x(\alpha)_{-1}, x(\alpha)_{0}$, and $x(\alpha)_{1}$, ($\alpha\in (\Bbb{C}^{3}\setminus 0)\cup\{\infty\}$) (cf. [Po, 3]).  It is clear that the set of subalgebras $S_{\alpha}$ is parameterized by $\Bbb{CP}^{2}\cup\{\infty\}$.  From the definition, it is also clear that $S_{\alpha}$ ($\alpha\in \Bbb{CP}^{2}\cup\{\infty\}$) are  right coideals of $\Cal{A}(\sigma)$, thus they are right quantum $\Cal{A}(\sigma)$-spaces, and we call them quantum super spheres.  The right coaction of $\Cal{A}(\sigma)$ on $\bold{x}(\alpha)$ is given by
$$
\Delta\bold{x}(\alpha)=\bold{x}(\alpha)\otimes M. \tag 4.3
$$
\par
The following lemma is originally due to Podle\`s [Po].
\proclaim{Lemma 4.2}
There are $\alpha$ with $\alpha_{0}\ne 0$ such that 
$$
\gather
c_{1}x(\alpha)_{-1}x(\alpha)_{1}+ c_{2}x(\alpha)_{1}x(\alpha)_{-1}
+ c_{3}x(\alpha)_{0}^{2} = b_{1}\cdot 1, \tag 4.4\\
c_{4}x(\alpha)_{-1}x(\alpha)_{1}+ c_{5}x(\alpha)_{1}x(\alpha)_{-1}
+ c_{6}x(\alpha)_{0}^{2} = b_{2}\sigma x(\alpha)_{0}+ b_{3}\cdot 1, \tag 4.5\\
m_{1}x(\alpha)_{-1}x(\alpha)_{0}+ m_{2}x(\alpha)_{0}x(\alpha)_{-1}
= m_{3}\sigma x(\alpha)_{-1}, \tag 4.6\\
n_{1}x(\alpha)_{1}x(\alpha)_{0}+ n_{2}x(\alpha)_{0}x(\alpha)_{1}
= n_{3}\sigma x(\alpha)_{1}, \tag 4.7
\endgather
$$
for nonzero $c_{i}, b_{j}, m_{k}, n_{k} \in \Bbb{C}$ ($1\le i\le 6$, $1\le j,k\le 3$).
\endproclaim
\demo{Proof}
The proof is given by solving the linear system obtained by writing (4.4)-(4.7) in the basis described by Theorem 1.1 and compare the coefficients on both the sides of the each of the equations. 
\enddemo
\par
Let us consider two special cases in more detail.  If $\alpha_{0}=0$, then
$$\gathered
x(\alpha)_{-1}= \alpha_{-1}a^{2} + \bold{i}\alpha_{1}c^{2},\\
 x(\alpha)_{0}= \alpha_{-1}\sqrt{1-q^{-1}}ab -\bold{i}\alpha_{1}\sqrt{1-q}dc,\\
x(\alpha)_{1}= \bold{i}\alpha_{-1}b^{2}+ \alpha_{1}d^{2}.
\endgathered \tag 4.8
$$
Thus, by using (1.1) we have
$$\split
x(\alpha)_{-1}x(\alpha)_{0} &= \alpha_{-1}^{2}\sqrt{1-q^{-1}}a^{3}b + \bold{i}\alpha_{-1}\alpha_{1}\sqrt{1-q}(t+t^{-2})abc^{2}\\
&\qquad + \bold{i}\alpha_{-1}\alpha_{1}\sqrt{1-q}ac\sigma - \alpha_{1}^{2}t^{-1}\sqrt{1-q}c^{3}d,
\endsplit
$$
and 
$$\split
x(\alpha)_{0}x(\alpha)_{-1} &= \alpha_{-1}^{2}t^{-2}\sqrt{1-q^{-1}}a^{3}b + \bold{i}\alpha_{-1}\alpha_{1}\sqrt{1-q}(1+t^{-5})abc^{2}\\
&\qquad + \bold{i}\alpha_{-1}\alpha_{1}t^{-2}\sqrt{1-q}ac\sigma - \alpha_{1}^{2}t^{-3}\sqrt{1-q}c^{3}d.
\endsplit
$$
Therefore there is no relation (4.6) possible in this case. Similarly, there is no relation (4.7) if $\alpha_{0}=0$.
\par
In the case $\alpha = \infty$, we have
$$\gather
x(\infty)_{0}^{2} - x(\infty)_{-1}x(\infty)_{1} + x(\infty)_{1}x(\infty)_{-1} = \sigma x(\infty)_{0}, \tag 4.9\\
x(\infty)_{0}^{2} + q^{-1}(1+q^{-1})x(\infty)_{-1}x(\infty)_{1} 
- (1+q^{-1})x(\infty)_{1}x(\infty)_{-1} = 1, \tag 4.10\\
qx(\infty)_{0}x(\infty)_{-1} - x(\infty)_{-1}x(\infty)_{0} = (1+q)\sigma x(\infty)_{-1}, \tag 4.11\\
x(\infty)_{0}x(\infty)_{1} - qx(\infty)_{1}x(\infty)_{0} = (1+q)\sigma x(\infty)_{1}. \tag 4.12
\endgather
$$
\par
Recall that if $A$ is a $\Bbb{Z}_{2}$-graded Hopf algebra, $(B, \psi)$ and $(B', \psi')$ are right quantum $A$-spaces, then we say that $B$ and $B'$ are isomorphic if there is an algebra isomorphism $\phi : B\longrightarrow  B'$ such that $\psi'\circ\phi = (\phi\otimes id)\circ \psi$.
\par
\proclaim{Theorem 4.3}(cf. [DK], [KS] and [Po]).
\par
(1) If $\bold{x}(\alpha)$ satisfies conditions (4.4)-(4.7), then the monomials
$$\gathered
x(\alpha)_{0}^{m}x(\alpha)_{-1}^{n}\sigma^{s}\quad (m,n \ge 0; s=0,1), \qquad\text{and} \\
x(\alpha)_{0}^{m}x(\alpha)_{1}^{n}\sigma^{s}\quad(m\ge 0, n\ge 1; s=0,1) 
\endgathered\tag 4.13
$$
form a basis of $S_{\alpha}$.
\par
(2) Two quantum super spheres $S_{\alpha}$ and $S_{\alpha'}$ are isomorphic if and only if there is an algebra homomorphism $\theta : S_{\alpha} \longrightarrow \Bbb{C}$ such that $\theta(\bold{x}(\alpha))M = \bold{x}(\alpha')C$ for some invertible $3\times 3$ complex matrix $C$.
\endproclaim
\demo{Proof}
The proof of (1) is similar to the non-super case, see for example [DK, p.308]. To prove (2), note that the condition for $S_{\alpha}$ and $S_{\alpha'}$ to be isomorphic is $\Delta\circ\phi = (\phi\otimes id)\circ\Delta$.  Suppose there exist a $\theta$ and an invertible matrix $C$ such that $\theta(\bold{x}(\alpha))M = \bold{x}(\alpha')C$.  We define an algebra homomorphism $\phi : S_{\alpha}\longrightarrow \Cal{A}(\sigma)$ by $\phi = m_{\Cal{A}(\sigma)}\circ (\theta\otimes id)\circ\Delta$.  Then by (4.3), we have
$$
\phi(\bold{x}(\alpha))=\theta(\bold{x}(\alpha))M = \bold{x}(\alpha')C.
$$
Since $\bold{x}(\alpha')C$ generates a quantum super sphere isomorphic to $S_{\alpha'}$, we conclude that $S_{\alpha}$ is isomorphic to $S_{\alpha'}$. Now suppose that $\phi : S_{\alpha}\longrightarrow S_{\alpha'}$ is an isomorphism. Computing different products of the generators of $S_{\alpha'}$ and expressed them in the basis of $\Cal{A}(\sigma)$ described in Theorem 1.1, it can be seen (compare with (4.4)-(4.7)) that since the presents of the element $\sigma$, there is an invertible matrix $C$ such that $\phi(\bold{x}(\alpha)) = \bold{x}(\alpha')C$.  Define $\theta : S_{\alpha}\longrightarrow \Bbb{C}$ by $\theta = \varepsilon\circ\phi$, then ($ m_{\Cal{A}(\sigma)}$ is omitted)
$$
\split
\theta(\bold{x}(\alpha))M &= (\theta\otimes id)
\circ\Delta(\bold{x}(\alpha))= (\varepsilon\otimes id )\circ (\phi\otimes id)\circ\Delta(\bold{x}(\alpha))\\
&= (\varepsilon\otimes id )\circ \Delta\circ\phi(\bold{x}(\alpha))
= (\varepsilon\otimes id )\circ\Delta(\bold{x}(\alpha')C)\\
&=\bold{x}(\alpha')C.
\endsplit
$$
\enddemo
\proclaim{Corollary 4.4}
The quantum super sphere $S_{\infty}$ is not isomorphic to any $S_{\alpha}$ with $\alpha \ne \infty$.
\endproclaim
\demo{Proof}
If $\theta : S_{\infty}\longrightarrow \Bbb{C}$ is an algebra homomorphism, let $y_{i}=\theta(x(\infty)_{i})$ ($i=-1,0,1$).  Since $\sigma x(\infty)_{0}= x(\infty)_{0}\sigma$ and $\sigma x(\infty)_{\pm 1}= -x(\infty)_{\pm 1}\sigma$, $(\sigma x(\infty)_{i})^{2}= x(\infty)_{0}^{2}$, hence $\theta(\sigma x(\infty)_{i})=\pm y_{i}$.  If $y_{-1}$ or $y_{1}$ is not equal to $0$, then (4.11) or (4.12) implies that $y_{0}= \pm \frac{1+q}{1-q}$, which contradicts (4.10).  Therefore $y_{-1}=y_{1}=0$ and $y_{0} = \pm 1$.  The corollary now follows from part (2) of Theorem 4.3.
\enddemo
\head
5. Finite dimensional representations
\endhead
In this and the following sections, we construct the finite dimensional comodules (left or right) of $\Cal{A}(\sigma)$ as bi-subcomodules of $\Cal{A}(\sigma)$, and express the matrix coefficients of these comodules in terms of the so-called little $q$-Jacobi polynomials. The computations involved in the discussion are quite similar to the non-super case, the main differences are due to the negative signs produced by the super commutativity and the formula $ad+tbc=\sigma$, we shall give details on where the differences occur and refer to the existing literature otherwise.
\par
Recall the following definition of the Gauss' binomial coefficients:
$$
(u;v)_{m}=\prod_{k=0}^{m-1}(1-uv^{k}), \quad 
\left(\matrix m\\ n \endmatrix\right)_{v}=\frac{(v;v)_{m}}{(v;v)_{n}(v;v)_{m-n}}. \tag 5.1
$$
For $m,n \in \Bbb{Z}_{+}$ we have
$$
\left(\matrix m+1\\ n+1 \endmatrix\right)_{v}=\left(\matrix m\\ n \endmatrix\right)_{v}v^{m-n} + \left(\matrix m\\ n+1 \endmatrix\right)_{v}, \tag 5.2
$$
and if $x$ and $y$ satisfy $xy=vyx$, we have
$$
(x+y)^{m}=\sum_{k=0}^{m}\left(\matrix m\\ k \endmatrix\right)_{v^{-1}}x^{k}y^{m-k}. \tag 5.3
$$
\par
By using (5.3) and induction on $m$ we can prove the following formulas
$$
\gathered
\Delta(a^{m}) =\sum_{k=0}^{m}(-1)^{[\frac{k}{2}]}\left(\matrix m\\ k \endmatrix\right)_{t^{-2}}a^{m-k}b^{k}\otimes a^{m-k}c^{k}, \\
\Delta(c^{m}) =\sum_{k=0}^{m}\left(\matrix m\\ k \endmatrix\right)_{t^{-2}}c^{m-k}d^{k}\otimes a^{m-k}c^{k},
\endgathered\tag 5.4
$$
where $[\frac{k}{2}]$ denotes the integer part of $\frac{k}{2}$, 
$$
\split
a^{m}d^{m} &=\sum_{k=0}^{m}\left(\matrix m\\ k \endmatrix\right)_{t^{-2}}
t^{2km-k^{2}}(cb)^{k}\sigma^{m-k}\\
&=\sum_{k=0}^{m}\left(\matrix m\\ k \endmatrix\right)_{t^{-2}}
t^{2km-k^{2}}(cb\sigma)^{k}\sigma^{m},
\endsplit\tag 5.5
$$
and
$$
\split
d^{m}a^{m}=\sum_{k=0}^{m}\left(\matrix m\\ k \endmatrix\right)_{t^{-2}}
t^{-k^{2}}(cb\sigma)^{k}\sigma^{m}.
\endsplit\tag 5.6
$$
\par
For each $\ell\in\frac{1}{2}\Bbb{Z}_{+}$, we let $I_{\ell} =\{-\ell,-\ell + 1, ...,\ell\}$ and define two $\Bbb C$-vector subspaces $V^{L}_{\ell}$ and $V^{R}_{\ell}$ of $\Cal{A}(\sigma)$ by (cf. [MNU, 2.1])
$$
V^{L}_{\ell}=\bigoplus_{i\in I_{\ell}}\Bbb{C}\xi^{(\ell)}_{i}\qquad\text{and}\qquad
V^{R}_{\ell}=\bigoplus_{i\in I_{\ell}}\Bbb{C}\eta^{(\ell)}_{i},\tag 5.7
$$
where (recall that $\bold{i}=\sqrt{-1}$)
$$
\xi^{(\ell)}_{i} =\bold{i}^{[\frac{\ell+i}{2}]}\left(\matrix 2\ell\\ \ell+i\endmatrix\right)^{\frac{1}{2}}_{t^{-2}}
a^{\ell -i}c^{\ell +i},\quad\text{and}\quad
\eta^{(\ell)}_{i} =\bold{i}^{[\frac{\ell+i}{2}]}\left(\matrix 2\ell\\ \ell+i\endmatrix\right)^{\frac{1}{2}}_{t^{-2}}
a^{\ell -i}b^{\ell +i}.\tag 5.8
$$
We shall denote the subspace spanned by the elements $\xi^{(\ell)}_{i}\sigma$ (resp. $\eta^{(\ell)}_{i}\sigma$) ($i\in I_{\ell}$) by $V^{L}_{\ell}\sigma$ (resp. $V^{R}_{\ell}\sigma$).
\par
As in the non-super case, we can prove that each $V^{L}_{\ell}\sigma^{s}$ (resp. $V^{R}_{\ell}\sigma^{s}$), $s=0,1$, forms an irreducible left (resp. right) $\Cal{A}(\sigma)$-subcomodule of $\Cal{A}(\sigma)$.  Let $m^{(\ell)}_{ij}\sigma^{s}$ ($i,j\in I_{\ell}$) be the matrix elements of $V^{L}_{\ell}\sigma^{s}$, $s=0,1$, then we have
$$
\Delta(\xi^{(\ell)}_{i}\sigma^{s})=\sum_{j\in I_{\ell}}m^{(\ell)}_{ij}\sigma^{s}\otimes\xi^{(\ell)}_{j}\sigma^{s}, \qquad i\in I_{\ell},\tag 5.9
$$
and
$$
\Delta(m^{(\ell)}_{ij}\sigma^{s})=\sum_{k\in I_{\ell}}m^{(\ell)}_{ik}\sigma^{s}\otimes m^{(\ell)}_{kj}\sigma^{s}, \qquad \text{and}\qquad
\varepsilon(m^{(\ell)}_{ij}\sigma^{s})=\delta_{ij}.\tag 5.10
$$
\par
The elements $\xi^{(\ell)}_{i}\sigma^{s}$ (resp. $\eta^{(\ell)}_{i}\sigma^{s}$) are the entries of the first column (resp. row) of the matrix $(m^{(\ell)}_{ij}\sigma^{s})$, (5.10) implies that (cf. Lemma 2.2 in [MNU])
$$
\Delta(\eta^{(\ell)}_{i}\sigma^{s})=\sum_{j\in I_{\ell}}\eta^{(\ell)}_{j}\sigma^{s}\otimes m^{(\ell)}_{ji}\sigma^{s}, \tag 5.11
$$
i.e. the right $\Cal{A}(\sigma)$-comodule $V^{R}_{\ell}\sigma^{s}$ is dual to the left $\Cal{A}(\sigma)$-comodule $V^{L}_{\ell}\sigma^{s}$. 
\par
For $\ell\in I_{\ell}$, $s=0,1$, let $M_{\ell}\sigma^{s}=\sum_{i,j\in I_{\ell}}\Bbb{C}m^{(\ell)}_{i,j}\sigma^{s}$, then $M_{\ell}\sigma^{s}$ is a bi-subcomodule of $\Cal{A}(\sigma)$, and we can define an $\Cal{A}(\sigma)$-bicomodule isomorphism $\Phi_{\ell} : V^{L}_{\ell}\sigma^{s}\otimes V^{R}_{\ell}\sigma^{s} \longrightarrow M_{\ell}\sigma^{s}$ by
$$
\Phi_{\ell}(\xi^{(\ell)}_{i}\sigma^{s}\otimes\eta^{(\ell)}_{j}\sigma^{s})= m^{(\ell)}_{ij}\sigma^{s}, \qquad
i,j\in I_{\ell}. \tag 5.12
$$
We also have $m^{(\ell)}_{ij}\sigma^{s}\in \Cal{A}(\sigma)[-2i,-2j]$, and for all $i,j\in \frac{1}{2}\Bbb{Z}$ such that $i-j \in \Bbb{Z}$,
$$
\Cal{A}(\sigma)[-2i,-2j]=\bigoplus_{\ell}(\Bbb{C}m^{(\ell)}_{ij}\oplus \Bbb{C}m^{(\ell)}_{ij}\sigma). \tag 5.13
$$
\par
We summarize our discussion in the following theorem, which can be proved as in [MNU, p. 367].
\proclaim{Theorem 5.1}
Any finite dimensional irreducible left (or right) $\Cal{A}(\sigma)$-comodule is isomorphic to $V^{L}_{\ell}\sigma^{s}$ (or $V^{R}_{\ell}\sigma^{s}$) ($s=0,1$) for some $\ell\in \frac{1}{2}\Bbb{Z}_{+}$, and the Hopf algebra $\Cal{A}(\sigma)$ is decomposed into the direct sum of $\Cal{A}(\sigma)$-bicomodules
$$
\Cal{A}(\sigma) =\bigoplus_{\ell\in\frac{1}{2}\Bbb{Z}_{+}}(M_{\ell}\oplus M_{\ell}\sigma).\tag 5.14
$$
\endproclaim
\par
\head 
6. Formulas for the matrix coefficients
\endhead
Since the matrix coefficients $m^{(\ell)}_{ij}\in\Cal{A}(\sigma)[-2i,-2j]$, by Proposition 1.3, there are $f^{(\ell)}_{ij}, g^{(\ell)}_{ij}\in\Cal{A}(\sigma)[0,0]$ such that
$$
m^{(\ell)}_{ij}=e_{-2i,-2j}f^{(\ell)}_{ij}=g^{(\ell)}_{ij}e_{-2i,-2j}. \tag 6.1
$$
To give explicit formulas for $m^{(\ell)}_{ij}$, we recall the definition of the little $q$-Jacobi polynomials (see [MNU, 2.2]):
$$
P^{(\alpha,\beta)}_{n}(z;q)=\sum_{r\ge 0}\frac{(q^{-n};q)_{r}(q^{\alpha+\beta+n+1};q)_{r}}{(q;q)_{r}(q^{\alpha+1};q)_{r}}(qz)^{r}, \tag 6.2
$$
and we use the following notation to simplify our statements (cf. [KS, 4.2.4]):
$$
N^{(\ell)}_{ij} = \bold{i}^{[\frac{\ell+i}{2}]-[\frac{\ell+j}{2}]} t^{(\ell+j)(i-j)}\left(\matrix \ell + i\\i-j\endmatrix\right)^{\frac{1}{2}}_{t^{-2}}\left(\matrix \ell -j\\i-j\endmatrix\right)^{\frac{1}{2}}_{t^{-2}}.\tag 6.3
$$
\par
\proclaim{Theorem 6.1}
For $\ell\in\frac{1}{2}\Bbb{Z}_{+}$, $i,j\in I_{\ell}$, the matrix coefficient $m^{(\ell)}_{ij}$ can be expressed in terms of the little $t$-Jacobi polynomials in $\zeta = tbc\sigma$ as follows.
\par
(1) If $i+j \le 0$, $i\ge j$, then 
$$
m^{(\ell)}_{ij}= a^{-i-j}c^{i-j}\sigma^{\ell+j}N^{(\ell)}_{ij}P^{(i-j,-i-j)}_{\ell+j}(\zeta; t^{-2}).\tag 6.4
$$
\par
(2) If $i+j \le 0$, $j\ge i$, then
$$
m^{(\ell)}_{ij}= (-1)^{[\frac{j-i}{2}]} 
a^{-i-j}b^{j-i}\sigma^{\ell+i}N^{(\ell)}_{ji}P^{(j-i,-i-j)}_{\ell+i}(\zeta; t^{-2}).\tag 6.5
$$
\par
(3) If $i+j \ge 0$, $j\ge i$, then
$$
m^{(\ell)}_{ij}= (-1)^{[\frac{j-i}{2}]} 
N^{(\ell)}_{-i,-j}P^{(j-i,i+j)}_{\ell-j}(\zeta; t^{-2})
b^{j-i}d^{i+j}\sigma^{\ell-j}.\tag 6.6
$$
\par
(4) If $i+j \ge 0$, $i\ge j$, then
$$
m^{(\ell)}_{ij}= 
N^{(\ell)}_{-j,-i}P^{(i-j,i+j)}_{\ell-i}(\zeta; t^{-2})
c^{i-j}d^{i+j}\sigma^{\ell-i}.\tag 6.7
$$
\endproclaim
\demo{Proof}
The proofs for cases (1) and (3) (resp. (2) and (4)) are somewhat symmetric to each other, we shall give proofs for (1) and (2).
\par
By (5.4) we have
$$
\split
\Delta(a^{\ell-i}c^{\ell+i}) &=\sum_{\overset {-\ell\le j\le\ell}\to{\mu\ge 0}}
\left(\matrix \ell-i\\ \mu \endmatrix\right)_{t^{-2}}
\left(\matrix \ell+i\\ \ell+j-\mu \endmatrix\right)_{t^{-2}}
(-1)^{[\frac{\mu}{2}]+\mu(i-j+\mu)}\\
&\times t^{-\mu(i-j+\mu)}a^{\ell-i-\mu}b^{\mu}c^{i-j+\mu}d^{\ell+j-\mu}
\otimes a^{\ell-j}c^{\ell+j}.
\endsplit
$$
If $i+j\le 0$ and $i\ge j$, by (5.5) the coefficient of $1\otimes a^{\ell-j}c^{\ell+j}$ is
$$
\split
&\sum_{\mu, k\ge 0}a^{-i-j}c^{i-j}
\left(\matrix \ell-i\\ \mu\endmatrix\right)_{t^{-2}}
\left(\matrix \ell+i\\ \ell+j- \mu\endmatrix\right)_{t^{-2}}
\left(\matrix \ell+j-\mu\\ k\endmatrix\right)_{t^{-2}}\\
&\qquad\times (-1)^{[\frac{\mu}{2}]+\mu+\frac{\mu(\mu-1)}{2}+k}
t^{\mu(\ell+j-\mu)+(i-j+\mu)(\ell+j-2\mu)+2k(\ell+j-\mu)-k^{2}}
(bc\sigma)^{k+\mu}\sigma^{\ell+j}\\
&\quad = t^{(\ell+j)(i-j)}a^{-i-j}c^{i-j}\sigma^{\ell+j}
\sum_{\mu, r\ge 0}
\left(\matrix \ell-i\\ \mu\endmatrix\right)_{t^{-2}}
\left(\matrix \ell+i\\ \ell+j- \mu\endmatrix\right)_{t^{-2}}\\
&\qquad\times \left(\matrix \ell+j-\mu\\ r-\mu\endmatrix\right)_{t^{-2}}t^{2r(\ell+j)-2\mu(i-j+\mu)-r^{2}}
(-bc\sigma)^{r},
\endsplit
$$
where at the last step we have set $r=\mu+k$ and used $[\frac{\mu}{2}]+\frac{\mu(\mu-1)}{2}\equiv 0(\mod 2)$.  Now by comparing with (2.27) and (2.29) in [MNU], the desired result follows in this case.
\par
If $i+j\le 0$ and $j\ge i$, we write the sum as
$$
\split
\Delta(a^{\ell-i}c^{\ell+i}) &=\sum_{\overset {-\ell\le j\le\ell}\to{\mu\ge 0}}
\left(\matrix \ell+i\\ \mu \endmatrix\right)_{t^{-2}}
\left(\matrix \ell-i\\ \ell-j-\mu \endmatrix\right)_{t^{-2}}
(-1)^{[\frac{j-i+\mu}{2}]+\mu(j-i+\mu)}\\
&\times t^{(j-i+2\mu)(\ell-i-\mu)-\mu(j-i+\mu)}
a^{-i-j}b^{j-i}b^{\mu}c^{\mu}a^{\ell+i-\mu}d^{\ell+i-\mu}
\otimes a^{\ell-j}c^{\ell+j}.
\endsplit
$$
The coefficient of $1\otimes a^{\ell-j}c^{\ell+j}$ is
$$
\split
&\sum_{\mu, k\ge 0}a^{-i-j}b^{j-i}\sigma^{\ell+i}
\left(\matrix \ell+i\\ \mu\endmatrix\right)_{t^{-2}}
\left(\matrix \ell-i\\ \ell-j- \mu\endmatrix\right)_{t^{-2}}
\left(\matrix \ell+i-\mu\\ k\endmatrix\right)_{t^{-2}}\\
&\qquad\times (-1)^{[\frac{j-i+\mu}{2}]+\mu(j-i+\mu)+\frac{\mu(\mu-1)}{2}+k}t^{(j-i+2\mu+2k)(\ell+i-\mu)-\mu(j-i+\mu)-k^{2}}
(bc\sigma)^{k+\mu}\\
&\quad = (-1)^{[\frac{j-i}{2}]}t^{(\ell+i)(j-i)}a^{-i-j}b^{j-i}\sigma^{\ell+i}
\sum_{\mu, r\ge 0}
\left(\matrix \ell+i\\ \mu\endmatrix\right)_{t^{-2}}
\left(\matrix \ell-i\\ \ell-j- \mu\endmatrix\right)_{t^{-2}}\\
&\qquad\times \left(\matrix \ell+i-\mu\\ r-\mu\endmatrix\right)_{t^{-2}}
t^{2r(\ell+i)-2\mu(j-i+\mu)-r^{2}}(-bc\sigma)^{r},\qquad\text{($r=\mu+k$)}
\endsplit
$$
where at the last step the sign was computed by
$$
\split
&[\frac{j-i+\mu}{2}]+\mu(j-i+\mu)+k+ \frac{\mu(\mu-1)}{2}\\
&\qquad \equiv\frac{(j-i)(j-i-1)}{2} + r\equiv[\frac{j-i}{2}] + r(\mod 2).
\endsplit
$$
Since
$$
\left(\matrix \ell+i\\ \mu\endmatrix\right)_{t^{-2}}
\left(\matrix \ell+i-\mu\\ r-\mu\endmatrix\right)_{t^{-2}}=
\left(\matrix \ell+i-\mu\\ r\endmatrix\right)_{t^{-2}}
\left(\matrix r\\ \mu\endmatrix\right)_{t^{-2}},
$$
for fixed $r$, we have
$$
\split
&\sum_{\mu\ge 0}
\left(\matrix \ell+i\\ \mu\endmatrix\right)_{t^{-2}}
\left(\matrix \ell-i\\ \ell-j- \mu\endmatrix\right)_{t^{-2}}
\left(\matrix \ell+i-\mu\\ r-\mu\endmatrix\right)_{t^{-2}}
t^{2r(\ell+i)-2\mu(j-i+\mu)-r^{2}}\\
&\qquad = \left(\matrix \ell-i+r\\ \ell-j\endmatrix\right)_{t^{-2}}
\left(\matrix \ell+i\\ r\endmatrix\right)_{t^{-2}}
t^{2r(\ell+i)-r^{2}}.
\endsplit
$$
Thus the coefficient of $1\otimes a^{\ell-j}c^{\ell+j}$ is
$$
(-1)^{[\frac{j-i}{2}]}t^{(\ell+i)(j-i)}a^{-i-j}b^{j-i}\sigma^{\ell+i}
\left(\matrix \ell-i\\ j-i\endmatrix\right)_{t^{-2}}
P^{j-i,-i-j}_{\ell+i}(tbc\sigma;t^{-2}),
$$
and after normalization, we get the formula in the theorem.
\enddemo
\par
\head
7. The Haar functional and the Peter-Weyl theorem
\endhead
We call a linear functional $\varphi : \Cal{A}(\sigma)\longrightarrow \Bbb{C}$ a left (resp. right) integral (cf. [MNU, 3.1] and [Mo, p.26]) if 
$$
(id\otimes \varphi)\circ\Delta(x)= \varphi(x)1_{\Cal{A}(\sigma)} \qquad
\text{(resp. $(\varphi\otimes id)\circ\Delta(x)= \varphi(x)1_{\Cal{A}(\sigma)}$)}\tag 7.1
$$
for any $x\in \Cal{A}(\sigma)$.
\par
Since the matrix elements $m^{(\ell)}_{ij}\sigma^{s}$ ($\ell \in \frac{1}{2}\Bbb{Z}_{+}, i,j\in I_{\ell}, s=0,1$) form a $\Bbb{C}$-basis of $\Cal{A}(\sigma)$, we can define a linear functional $h : \Cal{A}(\sigma)\longrightarrow \Bbb{C}$ by
$$
h(m^{(0)}_{00}\sigma^{s}) = h(\sigma^{s})=1, \qquad h(m^{(\ell)}_{ij})=0,\tag 7.2
$$
for all $\ell \in \frac{1}{2}\Bbb{Z}_{+}, i,j\in I_{\ell}, s=0,1$.  It follows from (5.10) that $h$ is both a left and a right integral on $\Cal{A}(\sigma)$, we call $h$ a Haar functional on $\Cal{A}(\sigma)$. The Haar functional $h$ can also be described in the basis $a^{i}b^{j}c^{k}d^{l}\sigma^{s}$ of $\Cal{A}(\sigma)$ provided by Theorem 1.1.  In fact, Theorem 3.1 of [MNU] holds in our case.
\proclaim{Theorem 7.1 ([KS, p.113], [MNU, Theorem 3.1])} The functional $h$ is the unique left and right integral on $\Cal{A}(\sigma)$ such that that $h(\sigma^{s})=1$ ($s=0,1$). The functional $h$ vanishes on the space $\Cal{A}(\sigma)[m,n]$ if $(m,n)\ne (0,0)$ and on the basis elements $a^{i}b^{j}c^{k}\sigma^{s}$ and $b^{j}c^{k}d^{i}\sigma^{s}$ if $i\ne 0$ or if $i=0$ but $j\ne k$. On the subalgebra $\Cal{A}(\sigma)[0,0]=\Bbb{C}<bc,\sigma>=\Bbb{C}<\zeta, \sigma>$ (recall that $\zeta=tbc\sigma$), $h$ is determined by
$$
h(\sigma)=1\quad\text{and}\quad h(\zeta^{n})=\frac{1-t^{-2}}{1-t^{-2(n+1)}}
=\frac{q+1}{q+(-1)^{n}q^{-n}}.\tag 7.3
$$
Furthermore, $h(S(x))=h(x)$ and $h(x^{\ast})=\overline{h(x)}$ for all $x\in\Cal{A}(\sigma)$.
\endproclaim
\demo{Proof}
Note that the key to the proof provided in [MNU] is formula (3.4) of [MNU], and we just need to establish a similar formula for $\Cal{A}(\sigma)$.  By using induction on $n$, it is easy to show that for any $n\in \Bbb{Z}_{+}$
$$
a^{n}d^{n}=(\zeta;t^{2})_{n}\sigma^{n}\qquad\text{and}\qquad 
d^{n}a^{n}=(t^{-2}\zeta;t^{-2})_{n}\sigma^{n}. \tag 7.4
$$
Let $P : \Cal{A}(\sigma)\longrightarrow\Cal{A}(\sigma)[0,0]$ be the projection according to the decomposition of Proposition 1.3, then we have
$$
\split
(id\otimes P)\circ\Delta(\zeta^{n}) &= (id\otimes P)((-1)^{[\frac{n}{2}]}t^{n}\Delta(b^{n})\Delta(c^{n})\Delta(\sigma^{n}))\\
&=(-1)^{[\frac{n}{2}]}t^{n}(id\otimes P)
((\sum_{i=0}^{n}(-1)^{i(n-i)}\binom{n}{i}_{t^{-2}}a^{i}b^{n-i}\otimes b^{i}d^{n-i})\\
&\qquad\times(\sum_{j=0}^{n}\binom{n}{j}_{t^{-2}}c^{j}d^{n-j}\otimes a^{j}c^{n-j})(\sigma^{n}\otimes\sigma^{n}))\\
&= (-1)^{[\frac{n}{2}]}t^{n}\sum_{j=0}^{n}\binom{n}{j}^{2}_{t^{-2}}
a^{j}b^{n-j}c^{n-j}d^{j}\sigma^{n}\otimes b^{j}d^{n-j}a^{n-j}c^{j}\sigma^{n}\\
&= (-1)^{[\frac{n}{2}]}t^{n}\sum_{j=0}^{n}\binom{n}{j}^{2}_{t^{-2}}t^{2j(n-j)}
(-1)^{\frac{(n-j)(n-j-1)}{2}+j(n-j)+\frac{j(j-1)}{2}}\\
&\qquad\times (bc\sigma)^{n-j}(\zeta;t^{2})_{j}\otimes (bc\sigma)^{j}(t^{-2}\zeta;t^{-2})_{n-j}.
\endsplit
$$
Since $\frac{(n-j)(n-j-1)}{2}+j(n-j)+\frac{j(j-1)}{2}\equiv [\frac{n}{2}](\mod 2)$, we have
$$
(id\otimes P)\circ\Delta(\zeta^{n})= \sum_{j=0}^{n}\binom{n}{j}^{2}_{t^{-2}}t^{2j(n-j)}(\zeta)^{n-j}(\zeta;t^{2})_{j}\otimes (\zeta)^{j}(t^{-2}\zeta;t^{-2})_{n-j}.\tag 7.5
$$
Now use (7.5) and argue as in the proof of Theorem 3.1 of [MNU].
\enddemo
\par
By the dual Maschke theorem (see [Mo, pp.25-26]) we have
\proclaim{Corollary 7.2}
The algebra $\Cal{A}(\sigma)$ is cosemisimple.
\endproclaim
\par
From now on to the end of this section we assume that $q<0$ is a real number. Then $t$ is also a real number. Recall that $\ast$ is defined by (1.6).  We use the Haar functional $h$ to define two hermitian forms $<,>_{R}$ and $<,>_{L}$ on $\Cal{A}(\sigma)$ by (cf. [MNU, (3.10)])
$$
<x,y>_{R}=h(xy^{\ast}), \qquad <x,y>_{L}=h(x^{\ast}y),\qquad x,y\in\Cal{A}(\sigma).\tag 7.6
$$
It follows from the definition that $<,>_{R}$ is conjugate linear in the second variable and $<,>_{L}$ in the first, and 
$$
<xy,z>_{R}=<x,zy^{\ast}>_{R}, \qquad <xy,z>_{L}=<y,x^{\ast}z>_{L}, \qquad x,y,z\in\Cal{A}(\sigma).\tag 7.7
$$
The fact that any finite dimensional left (resp. right) $\Cal{A}(\sigma)$-subcomodule of $\Cal{A}(\sigma)$ is unitary with respect to $<,>_{R}$ (resp. $<,>_{L}$) also follows from the definition of these two hermitian forms since $h$ is both a left and a right integral (cf. Proposition 3.5 of [MNU]). 
\par
To compute the values of these two hermitian forms on specific elements of $\Cal{A}(\sigma)$, recall the definition of $e_{mn}$ from (1.12).  We have (cf. [MNU, 3.17] and [KS, p.116])
$$
\gathered
e_{mn}e_{mn}^{\ast}= t^{\frac{(n-m)(n+m-2)}{2}}\zeta^{\frac{n-m}{2}}(\zeta;t^{2})_{\frac{m+n}{2}},\qquad \text{if $m+n\ge 0, m\le n$},\\
e_{mn}e_{mn}^{\ast}=t^{\frac{(m-n)(n+m-2)}{2}}\zeta^{\frac{m-n}{2}}(\zeta;t^{2})_{\frac{m+n}{2}},\qquad \text{if $m+n\ge 0, m\ge n$},\\
e_{mn}e_{mn}^{\ast}=\zeta^{\frac{m-n}{2}}(t^{-2}\zeta;t^{-2})_{\frac{-m-n}{2}},\qquad \text{if $m+n\le 0, m\ge n$},\\
e_{mn}e_{mn}^{\ast}=t^{m-n}\zeta^{\frac{n-m}{2}}(t^{-2}\zeta;t^{-2})_{\frac{-m-n}{2}},\qquad \text{if $m+n\le 0, m\le n$}.
\endgathered\tag 7.8
$$
These formulas can be derived from (1.12), (5.5) and (5.6).  For example, if $m+n\ge 0, m\le n$, we have
$$
\split
e_{mn}e_{mn}^{\ast}&=a^{\frac{m+n}{2}}c^{\frac{n-m}{2}}(-t^{-1}
b\sigma)^{\frac{n-m}{2}}(d\sigma)^{\frac{m+n}{2}}= a^{\frac{m+n}{2}}(t^{-1}
bc\sigma)^{\frac{n-m}{2}}(d\sigma)^{\frac{m+n}{2}}\\
&= t^{\frac{(n-m)(n+m-2)}{2}}\zeta^{\frac{n-m}{2}}(\zeta;t^{2})_{\frac{m+n}{2}}.
\endsplit
$$
\par
We rewrite formulas (3.26) and (3.27) of [MNU] (cf. [KS, p. 112]) in our case as 
$$
\gathered
h(\zeta^{r}(\zeta;t^{2})_{s}) = t^{-2(r+1)}\frac{(t^{-2};t^{-2})_{r}(t^{-2};t^{-2})_{s}(t^{-2};t^{-2})_{1}}{(t^{-2};t^{-2})_{r+s+1}},\\
h(\zeta^{r}(t^{-2}\zeta;t^{-2})_{s}) = \frac{(t^{-2};t^{-2})_{r}(t^{-2};t^{-2})_{s}(t^{-2};t^{-2})_{1}}{(t^{-2};t^{-2})_{r+s+1}},\\
\endgathered\tag 7.9
$$
\par
By using formulas (7.8) and (7.9), we can compute the values of $h$ and the two hermitian forms on the elements $\xi^{(\ell)}_{i}\sigma^{s}$, $\eta^{(\ell)}_{i}\sigma^{s}$ and $m^{(\ell)}_{ij}\sigma^{s}$. 
\par
\proclaim{Theorem 7.3} The algebra $\Cal{A}(\sigma)$ has an orthogonal decomposition
$$
\Cal{A}(\sigma) =\bigoplus_{\ell\in\frac{1}{2}\Bbb{Z}_{+}}(M_{\ell}\oplus M_{\ell}\sigma)\tag 7.10
$$
with respect to $<,>_{R}$ or $<,>_{L}$. Furthermore, the following relations hold
$$\gather
<m^{(\ell)}_{ij}\sigma^{s},m^{(\ell')}_{i'j'}\sigma^{s'}>_{R} = [2\ell+1]^{-1}_{t}t^{2j}\delta_{\ell\ell'}\delta_{ii'}\delta_{jj'}, \tag 7.10\\
<m^{(\ell)}_{ij}\sigma^{s},m^{(\ell')}_{i'j'}\sigma^{s'}>_{L} = [2\ell+1]^{-1}_{t}t^{-2i}\delta_{\ell\ell'}\delta_{ii'}\delta_{jj'}, \tag 7.11
\endgather
$$
where $[2\ell+1]_{t}=\frac{t^{2\ell+1}-t^{-2\ell-1}}{t-t^{-1}}$.
\endproclaim
\demo{Proof}
Note that Lemma 3.8 of [MNU] holds in our case, so we can apply the argument in the proof of Theorem 3.7 of [MUN] to our case (cf. [KS, 4.3.2]).
\enddemo
\par

\enddocument